\documentclass{amsart}

\title[Attracting Basins]{Attracting Basins of Volume
Preserving Automorphisms of $\CC^k$}
\author{Han Peters, Liz Raquel Vivas, Erlend Forn\ae ss Wold}
\date{August 20, 2006}

\usepackage{amscd}

\newtheorem{theorem}{Theorem}
\newtheorem{lemma}{Lemma}
\newtheorem{proposition}{Proposition}
\newtheorem{corollary}{Corollary}
\newtheorem{example}{Example}

\theoremstyle{definition}

\theoremstyle{remark}
\newtheorem{remark}{Remark}

\newcommand{\NN}{\mathbb{N}}
\newcommand{\RR}{\mathbb{R}}
\newcommand{\CC}{\mathbb{C}}

\newcommand{\autidentity}{\mathrm{Aut}_1^1(\mathbb{C}^2, 0)}

\def\a{{\alpha}}

\def\d{{\delta}}
\def\e{{\epsilon}}
\def\s{{\sigma}}

\def\r{{\rho}}
\def\l{{\lambda}}
\def\O{{\Omega}}
\def\Aut{{\mathrm{Aut}}}

\begin{document}

\bibliographystyle{plain}

\begin{abstract}
We study topological properties of attracting sets for
automorphisms of $\CC^k$. Our main result is that a generic volume
preserving automorphism has a hyperbolic fixed point with a dense
stable manifold. We prove the same result for volume preserving
maps tangent to the identity. On the other hand, we show that an
attracting set can only contain a neighborhood of the fixed point
if it is an attracting fixed point. We will see that the latter
does not hold in the non-autonomous setting.
\end{abstract}

\maketitle

\section{introduction}

Let $f$ be an automorphism of $\CC^k$ with a fixed point at the
origin. Even if the origin is not an attracting fixed point, there
can still be points whose orbits converge to the origin. In this
paper we will study how large such an attracting set can be. We
will make this more precise later.

The behavior of an attracting set varies greatly depending on the
eigenvalues of $\mathrm{d}f(0)$. If all eigenvalues have modulus
strictly smaller than $1$ then we say that $f$ has an
{\emph{attracting}} fixed point. This is the easiest situation,
the attracting set must contain a neighborhood of the origin and
is biholomorphic to $\CC^k$ \cite{rr}. The situation is similar
when all eigenvalues have modulus strictly larger than $1$, one
just considers the inverse mapping.

The fixed point is called {\emph{hyperbolic}} if no eigenvalues
have modulus $1$, and there are eigenvalues of modulus greater
than $1$ as well as less than $1$. In this case the attracting set
is biholomorphic to $\CC^m$, where $m$ is the number of
eigenvalues of modulus less than $1$ (this follows from
\cite{rr}).

The complex structure of attracting sets has also been studied in
the {\emph{semi-attracting}} case (eigenvalues of modulus smaller
than and equal to $1$) , and for automorphisms {\emph{tangent to
the identity}} (where $\mathrm{d}f(0) = \mathrm{Id}$). In both
cases the attracting set can also be biholomorphically equivalent
to (possibly lower dimensional) complex Euclidean space, see for
example \cite{ue} for the semi-attracting case and \cite{we},
\cite{hak1} and \cite{hak2} for automorphisms tangent to the
identity.

In this article we do not study the complex structure of
attracting sets but instead we look at topological properties.
Suppose an automorphism has a fixed point that is not attracting
but does have a non-trivial attracting set $\O$. We are interested
in three related questions:

\noindent (a) Can $\O$ be dense?

\noindent (b) Can $\O$ have interior points?

\noindent (c) Can $\O$ contain a neighborhood of the origin?

Our main result is an affirmative answer to Question (a). More
precisely, we will show the following

\begin{theorem}\label{dense}
There is a dense $G_\d$-set $\mathcal{V}$ of volume preserving
automorphisms of $\CC^k$ that have a hyperbolic fixed point whose
stable manifold is dense in $\CC^k$.
\end{theorem}

Here $\mathcal{V}$ is a dense $G_\d$-subset of the set of volume
preserving automorphisms of $\CC^k$, equipped with the compact
open topology.

We will then focus on volume preserving automorphisms that are
tangent to the identity. We will prove the existence of dense
attracting sets for these maps along the same lines as for a
hyperbolic fixed point.

The answer to Question (b) is obvious, since it is possible to
have an attracting set that is biholomorphic to $\CC^k$, where $k$
is the dimension of the ambient space \cite{ue}, \cite{we},
\cite{hak2}. However, we will easily see that the attracting set
of a {\emph{volume preserving}} automorphism cannot have interior
points.

We will also show that the answer to Question (c) is negative, if
the attracting set contains a neighborhood of the origin then the
fixed point must be attracting. We note that this result depends
on the holomorphicity of the mapping, as well as on the ambient
space $\CC^k$.

Finally, we will see that in the non-autonomous setting the basin
can be all of $\CC^k$, even if all the mappings are tangent to the
identity.

In the next section we will set notation and answer question (c).
We will prove our main result in Section 3, and show the analogous
statement for maps tangent to the identity in Section 4. In the
last section we will treat the non-autonomous setting.

\section{Fixed Point in the Interior}

We denote by $\Aut(\CC^k)$ the set of holomorphic automorphisms of
$\CC^k$ and by $\Aut_1(\CC^k)$ the set of volume preserving
automorphisms of $\CC^k$, both equipped with the compact-open
topology.

We let $\|\cdot\|$ denote the Euclidean norm on $\CC^k$ and for
$r>0$ we write $B_r(z) \subset \CC^k$ for the ball of radius $r$
centered at $z$. When $z = 0$ we will just write $B_r$.

For $f \in \Aut(\CC^k)$ with a fixed point $p$ we will study the
{\emph{attracting set}}

$$
\O = \{z \in \CC^k \mid f^n(z) \rightarrow p, \; \mathrm{as} \; n
\rightarrow \infty\}
$$

When $p$ is an attracting fixed point this attracting set is
generally called the attracting basin, and when $p$ is a
hyperbolic fixed point it is called the stable manifold.

The automorphisms of $\CC^k$ constructed by Ueda, Hakim and
Weickert that have a neutral or semi-attractive fixed point with
an attracting set biholomorphic to $\CC^k$ all have the fixed
point lying in the boundary of the basin. It is natural to ask
whether the attracting set of such a fixed point can ever contain
an open neighborhood of the fixed point. The following result
shows that this cannot happen.

\begin{theorem}\label{uniform}
Let $f:\CC^k\rightarrow\CC^k$ be a holomorphic map such that
$f(0)=0$, and let $\O$  be the attracting set. If $\O$ contains a
neighborhood of the origin then $0$ is an attracting fixed point.
\end{theorem}
\begin{proof}
Suppose for the purpose of a contradiction that the closed ball
$\overline{B}_r$ is in the attracting set for some $r>0$, and that
$0$ is not an attracting fixed point. It follows that no iterate
of $f$ has $0$ as an attracting fixed point.

Our first claim is that the set
$\{f^n(\overline{B}_r)\}_{n\in\NN}$ is unbounded. If not then the
set of iterates $\{f^n\}$  would be a normal family on
$\overline{B}_r$, and we may pass to a convergent subsequence. But
since $f^n(x)\rightarrow\ 0$ for all $x \in \overline{B}_r$ we get
that $f^N(\overline{B}_r)\subset\subset \overline{B}_r$  for some
$N\in\NN$, and by the contraction principle (see page 219 of
\cite{mo}) this contradicts the fact that $0$ is not an attracting
fixed point for $f^N$.

For $m\in\NN$ define the bounded set $K_m:=\cup_{i=1}^m
f^i(\overline{B}_r)$. It follows from the above claim that there
must be a point $x\in \overline{B}_r$ such that
$f^n(x)\in\CC^k\setminus K_m$ for some $n\in\NN$, so it follows
that in the sequence $\{f(x),f^2(x),...,f^n(x)\}$ there have to be
at least $m$ successive $f^i(x)$'s such that $f^i(x)\notin B_r$.
So each set $C_m=\{x\in \overline{B}_r \mid f^j(x)\notin B_r,
j=1,...,m\}$ is non-empty.  We have that
$C_1=f^{-1}(\CC^k\setminus B_r)\cap \overline{B}_r$, and then
$C_m=f^{-m}(\CC^k\setminus B_r)\cap C_{m-1}$ for $m=2,3,...$, so
we have a decreasing sequence of compact sets. Therefore there is
a point $x\in\cap_{i=1}^\infty C_i$, and it follows that $f^j(x)$
does not converge to the origin, which is a contradiction.
\end{proof}

\begin{example}
Theorem \ref{uniform} does not hold for holomorphic self maps of
complex manifolds in general.  The mapping $f(z)=\frac{z}{1+z}$ is
an automorphism of the Riemann sphere and the orbit $f^m(z)$ is
given by $f^m(z)=\frac{z}{1+mz}$.  So we see that the origin is an
attracting fixed point for $f$, and the attracting set is in fact
equal to the entire Riemann sphere. But since
$f^{m}(\frac{-1}{m})=\infty$  we see that the attraction is not
uniform in any neighborhood of the origin.
\end{example}

\begin{example}
Theorem \ref{uniform} obviously does not hold for diffeomorphisms,
consider for example (for $x \in \RR^k$)

\begin{align*}
x \mapsto (1-\|x\|^\frac{1}{\|x\|})x, \; \mathrm{for} \; x \neq 0,\\
0 \mapsto 0.
\end{align*}

Here the basin is the unit ball, and the attraction is uniform on
compact subsets. This raises the following question: If $f$ is a
homeomorphism of $\RR^n$, and suppose that $f$ has a fixed point
such that the attracting set contains a neighborhood of the fixed
point. Is the attraction necessarily uniform on compact subsets?
The answer to this question is also negative when $n \ge 2$.

For $x \in \RR$ let
$$
\psi(x) = \frac{x(4\pi - x)}{2\pi}.
$$
Notice that $\psi(0) = 0, \psi(2\pi) = 2\pi$ and $\psi^n(x)
\rightarrow 2\pi$ for any $x \in (0, 2\pi)$ as $n \rightarrow
\infty$.

Now let $f : \CC \rightarrow \CC$ be defined as follows:

For $r \ge 1$ and $\theta \in [0, 2 \pi)$, we define
$$
f(re^{i\theta}) = \frac{r+1}{2}e^{i\psi(\theta)}.
$$

To define $f(z)$ for $z$ inside the unit ball, note that any such
$z$ lies on a unique circle through $1$ that is tangent to the
unit circle. Let $f$ fix those circles, so that the angle (with
respect to the center of such a circle) of $f(z)$ becomes
$\psi(\theta)$, where $\theta$ is the angle of $z$.

Note that $f$ is continuous, has a fixed point at $1$ and $f^n(z)$
converges to $1$ for any $z \in \CC$. Yet the convergence is not
uniform. With a little care the same construction works for a
diffeomorphism.
\end{example}

\section{Dense Stable Manifolds}

Before we prove Theorem \ref{dense}, we first show that an
attracting set of a volume preserving automorphism cannot have
interior points. We have already noted that generally the
attracting set of a semi-attracting or neutral fixed point can be
biholomorphic to $\CC^k$, so can in particular have interior.
However, we easily see that this cannot be the case when dealing
with volume preserving automorphisms.

\begin{proposition}\label{interior}
Let $f \in \Aut_1(\CC^k)$ have a fixed point at the origin, and
let $\Omega$ be the attracting set. Then $\Omega$ has empty
interior.
\end{proposition}
\begin{proof}
Let $\e > 0$ and define $U$ as the set of those points whose
forward orbit lies entirely in $B_\e$. Then $U$ is forward
invariant under $f$, and for every $n \in\NN$ we have that $f^n
(U) \subset B_\e$ and $\mathrm{Vol}(f^n(U)) = \mathrm{Vol}(U)$. It
follows that for every $n \in \NN$ we have that the set $\{z
\notin B_\e, f^n(z) \in U\}$ has no volume. So the countable union
$$
\{ z \notin B_\e \mid \exists n \in \NN: f^n(z) \in U\}
$$
has empty interior. But the orbit of any point in $\Omega$ must
eventually land in $U$. Hence the set $\Omega \setminus B_\e$ has
empty interior. Since this holds for any $\epsilon >0$ the proof
is complete.
\end{proof}

For an automorphism $f$ with a hyperbolic fixed point $p$ denote
the {\emph{local stable manifold}} by
$$
\Sigma^f_\epsilon(p) = \{ z \in \CC^k \mid \|f^n(z) - p\| < \e \
\forall n \in \NN \}.
$$
For small enough $\e$ we have that $\Sigma^f_\e(p)$ is a graph
over the attracting direction of $\mathrm{d}f(p)$, and
$\{f^n(z)\}$ converges to $p$ if and only if $f^n(z) \in
\Sigma^f_\e(p)$ for some $n \in \NN$ (see for instance Chapter 6.2
in Katok-Hasselblat \cite{kh}). In other words, if we denote the
attracting set or {\emph{stable manifold}} by $\Sigma^f(p)$ then
$$
\Sigma^f(p) = \bigcup_{n\in\NN} f^{-n} \Sigma^f_\e(p).
$$
As noted before, it follows from the appendix of \cite{rr} that
$\Sigma^f(p)$  is biholomorphic to $\CC^m$, where $m$ is the
number of attracting directions.

To prove Theorem \ref{dense} we need a stability condition for
stable manifolds.  Let us fix an $f\in \mathrm{Aut}_1(\CC^k)$ with a
fixed point at the origin, and assume that $f$ is of the form
$$
f(z)=(\l_1 z_1 + \a_1(z), \ldots ,\l_m z_m + \a_m(z),\mu_1 z_{m+1}
+ \a_{m+1}(z),\ldots,\mu_{k-m} z_k +\a_k(z)),
$$
where $|\l_i|<1$, $|\mu_i|>1$, and the $\a_i$'s are functions of
degree at least two. For $\d>0$ we let $\triangle^m_\d$ denote the
polydisk $\triangle^m_\d=\{(z_1,\cdots,z_k)\in\CC^k \mid z_i=0 \
\mathrm{for} \ i>m,|z_i|<\d \ \mathrm{for} \ i\leq m\}$. As stated
above: If $\e$ is small enough then for all $\d<\e$  we have that
$\Sigma_\e^f(0)$ is (locally) a graph $\Gamma^f_\d$ over
$\triangle^m_\d$.  We need the following proposition:

\begin{lemma}\label{graphs}
Let $\{f_j\}\subset \mathrm{Aut}_1(\CC^k)$  such that
$\|f-f_j\|_{\overline\triangle_\e^k}\rightarrow 0$  as
$j\rightarrow\infty$. Then there exists a fixed $\d<\e$ such that
for all $j$ large enough:

\

\noindent (a) $f_j$  has a unique hyperbolic fixed point $p_j$  in
$\triangle^k_\e$, and $p_j\rightarrow 0$  as $j\rightarrow\infty$,
\

\noindent (b) $\Sigma_\e^{f_j}(p_j)$  is (locally) a graph
$\Gamma_\d^{f_j}$ over $\triangle^m_\d$, \

\noindent (c)
$\mathrm{d}(\Gamma_\d^{f_j},\Gamma_\d^{f})\rightarrow 0$ as
$j\rightarrow\infty$, where $\mathrm{d}(\cdot,\cdot)$ denotes the
Hausdorff distance.
\end{lemma}

\emph{Sketch of the proof.} (a) is well known, and for (b) and (c)
we may well assume that $p_j = o$ for high enough $j$.

In a small enough polydisc $\Delta_\d$, the map $f$ is strictly
expanding in the repelling directions and strictly contracting in
the attracting direction. This gives that $\Sigma_\d^{f}(0)$ is a
graph over the attracting direction. For $f_j$ close enough to $f$
we have that $f_j$ is still strictly expanding and contracting in
this same polydisc and we get (b).

For $\gamma >0$ arbitrarily small, let $\mathcal{N}_\gamma$ be the
$\gamma$ neighborhood of $\Sigma_\d^{f}(0)$, and let $K =
\Delta_\d - \mathcal{N}_\gamma$. Then there is an $n \in \NN$ such
that for every $z \in K$ there is an $j \in \{1, \cdots , n\}$
such that $f^j(z) \notin \Delta_\d$. Hence the same is true for
$f_j$ close enough, so we have that $\Sigma_\d^{f_j}(0) \subset
\mathcal{N}_\gamma$. But since $\Sigma_\d^{f_j}(0)$ and
$\Sigma_\d^{f}(0)$ are both graphs over the stable direction, we
must have that $\mathrm{d}(\Sigma_\d^{f}(0), \Sigma_\d^{f_j}(0)) <
\gamma$. \hfill $\Box$

\

For each $n\in\NN$  we now let $\Gamma_\d^{f_j}(n)$  denote the
set $f_j^{-n}(\Gamma_\d^{f_j})$, such that
$\Sigma^{f_j}(p_j)=\bigcup_{n\in\NN}\Gamma_\d^{f_j}(n)$. The
following is then an immediate consequence of the above proposition:

\begin{corollary}\label{close}
Let $U$  be any neighborhood of
$\overline\triangle_\e^k\cup\Gamma_\d^f(n)$  for some fixed
$n\in\NN$, and let $\{f_j\}\subset \mathrm{Aut}_1(\CC^k)$  such
that $\|f-f_j\|_{\overline U}\rightarrow 0$  as
$j\rightarrow\infty$. Then
$\mathrm{d}(\Gamma_\d^{f_j}(n),\Gamma_\d^f(n))\rightarrow 0$ as
$j\rightarrow\infty$.
\end{corollary}

\begin{proposition}\label{perturb}
Let $f\in \mathrm{Aut}_1(\CC^k)$  as above have a hyperbolic fixed
point at the origin, let $\d,\r>0$, let $q\in\CC^k$ and let $K$ be
a compact subset of $\CC^k$.  Then there exists a $g\in
\mathrm{Aut}_1(\CC^k)$  such that the following hold:

\

(a) $g$  has a unique hyperbolic fixed point $p$  close to the
origin, \

(b) $\|g-f\|_K<\d$, \

(c) There is a point $q'\in\Sigma^g(p)$  such that $\|q'-q\|<\r$.
\end{proposition}
\begin{proof}
We assume that (c) is not already satisfied by $\Sigma^f(0)$, and
we assume that $K$  is a closed ball.  By Theorem 3.1 in
\cite{fo2} there is a $\tilde g\in \mathrm{Aut}_1(\CC^k)$ with
$\|\tilde g-f\|_K<\frac{\d}{2}$  and such that the unbounded
orbits of $\tilde g$  are dense in $\CC^k$. Choose $\tilde
q\in\CC^k$  with $\|\tilde q-q\|<\frac{\r}{2}$ such that $\{\tilde
g^n(\tilde q)\}_{n\in\NN}$ is unbounded. By Lemma \ref{graphs} we
know that if $\tilde{g}$ is a good enough approximation of $f$
then $g$ has a hyperbolic fixed point $\tilde p$ near $0$. We have
that $\Sigma^{\tilde g}(\tilde p)$ is unbounded since
$\Sigma^{\tilde g}(\tilde p)$ is biholomorphic to $\CC^m$ (This
follows from the appendix of \cite{rr}).

So we may choose a point $x\in\Sigma^{\tilde g}(\tilde{p})$  such
that $x\in\CC^k\setminus K$  and such that $\tilde g^n(x)\in K$
for all $n\geq 1$.  Let $M$ be an integer such that $\tilde
g^M(x)\in\Gamma_\d^{\tilde g}$, and let $N$ be the smallest
integer such that $\tilde g^N(\tilde q)\subset\CC^k\setminus K$.
For some $r>0$ we have that the set $K\cup B_r(x)\cup B_r(\tilde
g^N(\tilde q))$ is polynomially convex, and we let $\phi\in
\mathrm{Aut}_1(\CC^k)$  such that $\phi(B_r(\tilde g^N(\tilde
q))=B_r(x)$. Let $\mathcal{N}$ be a small enough neighborhood of
$\tilde q$ such that $\tilde g^N(\mathcal N)\subset\subset
B_r(\tilde g^N(\tilde q))$, and let $\mathcal V$ be a small enough
neighborhood of $\tilde g^M(x)$ such that $\mathcal
V\subset\subset\tilde g^M\circ\phi\circ\tilde g^N(\mathcal N)$. \

By \cite{fr} and \cite{fr2} there is exists a sequence of
automorphisms $\phi_j\in \mathrm{Aut}_1(\CC^k)$  such that
$\phi_j\rightarrow\phi$  on $B_r(\tilde g^N(x))$  and such that
$\phi_j\rightarrow\mathrm{Id}$ on $K$.  Approximating by a
\emph{volume preserving} automorphism is possible because of the
vanishing of the following cohomology group \cite{fr2}
$$
H^{k-1}(K\cup\overline{B_r(\tilde g^N(\tilde q))},\CC)=0.
$$

Let $\Phi_j$  denote $\tilde g\circ\phi_j$.  If $j$  is large
enough we have that $\mathcal
V\subset\subset\Phi_j^{M+N+1}(\mathcal N)$, and by Lemma
\ref{graphs} we have that $\Gamma_{\d}^{\Phi_j}\cap\mathcal
V\neq\emptyset$  if $j$  is large.  So the global stable manifold
of $\Phi_j$ intersects $\mathcal N$, and the result follows by
letting $g=\Phi_j$.
\end{proof}

\begin{corollary}\label{open}
Let $q_1,...,q_m$  be points in $\CC^k$  and let $\e>0$.  Then the
set of automorphisms $f\in \mathrm{Aut}_1(\CC^k)$  having a stable
manifold $\Sigma^f(p)$  with points $p_1,...,p_m\in\Sigma_p^f$ such
that $\|p_i-q_i\|<\e$  is dense and open.
\end{corollary}
\begin{proof}
Let $f\in \mathrm{Aut}_1(\CC^k)$.  Let $N\in\NN$  be arbitrary,
$\r>0$, and choose $p\notin(B_N\cup f(B_N))$.  Since
$H^{k-1}(f(B_N),\CC)=0$ there is a sequence of automorphisms
$\{g_j\} \in \mathrm{Aut}_1(\CC^k)$  with $g_j(f(p))=p$  and such
that $\|g_j-\mathrm{Id}\|_{f(B_N)}\rightarrow 0$ \cite{fr}.  By
composing with a linear map arbitrarily close to the identity if
necessary, we may assume that each $g_j\circ f$ has a hyperbolic
fixed point at $p$.  If $j$  is large we have that $\|g_j\circ f -
f\|_{B_N}<\r$, and it follows that the set of volume preserving
automorphisms having a hyperbolic fixed point is dense. By
Proposition \ref{dense}  it is also open.  \

Note that by Corollary \ref{close} the set of volume preserving
automorphisms having a stable manifold with a point $p_i$ that is
$\e$-close to some point $q_i$  is open.  Therefore it is enough
to consider the point $q_1$.  Let $h$  denote $g_j\circ f$ for a
some large $j$.  By Proposition \ref{perturb}  there exists for
any $\r>0$  a $\tilde h\in\mathrm{Aut}_1(\CC^k)$  such that
$\|\tilde h-h\|_{B_N}<\r$  and such that $\|p_1-q_1\|<\r$  for
some $p_1\in\Sigma^{\tilde{h}}_p$, where $p$  is a hyperbolic
fixed point for $\tilde h$.  The result follows.

\end{proof}

In the following proof, note that $\mathrm{Aut}_1(\CC^k)$  is a
\emph{Baire Space}, meaning that a countable intersection of open
and dense sets is again dense.

\

\emph{Proof of Theorem \ref{dense}.} Let $\{q_i\}_{i\in\NN}$ be a
dense set of points in $\CC^k$ and let $\e_j\searrow 0$.  For each
$j\in\NN$ let $\mathcal{V}_j$ denote the set of automorphisms
$f\in \mathrm{Aut}_1(\CC^k)$  such that $f$ has a stable manifold
$\Sigma_p^f$ with points $p_1,...,p_j\in\Sigma_p^f$  and
$\|p_i-q_i\|<\e_j$. According to Corollary \ref{open} each
$\mathcal{V}_j$ is open and dense. Since
$\mathcal{V}:=\cap_{j\in\NN}\mathcal{V}_j$ is dense the result
follows. \hfill $\Box$

\section{Automorphisms of $\CC^2 $ tangent to the identity}

We will now show that dense attracting sets also occur for volume
preserving automorphisms that are tangent to the identity. We will
restrict ourselves to automorphisms of $\CC^2$, and we will see
that a statement analogous to Theorem \ref{dense} holds for volume
preserving automorphisms tangent to the identity. Since the proof
is almost identical to the proof of Theorem \ref{dense} we will
not show it in great detail. The main difficulty is to prove that
the attracting set is (locally) stable under small perturbations.

We let $\autidentity$ be the set of volume preserving
automorphisms of $\CC^2$ that are tangent to the identity. We
equip $\autidentity$ with the compact open topology as before. For
$f \in \autidentity$ we write $f(z) = z + P_2(z) + \ldots$ where
$P_2(z)$ is homogeneous of degree $2$. Recall from \cite{hak1}
that $v \in \CC^2$ is called a characteristic direction if $P_2(v)
= \lambda v$ for some $\lambda \in \CC$. If $\lambda \neq 0$ then
$v$ is a \emph{non-degenerate} characteristic direction.

We first claim that every automorphism tangent to the identity
must have a characteristic direction. If $P_2(z) \equiv 0$ then it
is clear, so without loss of generality we assume that $P_2 =
(p,q)$ with $p(x,y) \neq 0$ for some $(x,y)$. We blow up $y = ux$
and we get that $P_2$ gives the rational function
\begin{align}\label{equation}
u \mapsto \frac{q(x,ux)}{p(x, ux)} = \frac{q(1,u)}{p(1, u)}.
\end{align}
Note that this function must have a fixed point, which is a
characteristic direction.

The subset of $\autidentity$ consisting of those automorphisms
that have a non-degenerate characteristic direction is open and
dense. More precisely, we have

\begin{lemma}
Let $f= \mathrm{Id} + P_2 + h.o.t. \in \autidentity$ with $P_2(v)
= \lambda v$. Then

\noindent (a) If $\lambda \neq 0$ and $\tilde{f}$ is close enough
to $f$ then $\tilde{f}$ has a has a non-degenerate characteristic
direction arbitrarily close to $v$.

\noindent (b) If $\lambda = 0$ then there exist mappings in
$\autidentity$ arbitrarily close to $f$ that have $v$ as a
non-degenerate characteristic direction.
\end{lemma}
\begin{proof}
Assertion (a) follows immediately from the fact that the
characteristic directions are fixed points of the rational
equation \eqref{equation}, whose parameters depend continuously on
the mapping $f$.

To prove assertion (b), assume without loss of generality that $v
= (1,0)$, in other words so that $P_2(1,0) = (0,0)$. Let $\psi_\e
= \mathrm{Id} + \Psi_\e + h.o.t. \in \autidentity$ with
$\Psi_\e(x,y) = (\e x^2m, -2\e xy)$, the existence of such an
automorphism follows from \cite{bf}. Then $\psi_\e \circ f$ has
$(1,0)$ as non-degenerate characteristic direction.
\end{proof}

Now let $f \in \autidentity$ have a non-degenerate characteristic
direction $v$. After a suitable conjugation by an affine mapping
we have that $v = (1,0)$ and that $P_2(v) = v$, i.e. $\lambda =
1$.

Write $P_2(z) = (p(z), q(z))$. Since $f$ is volume preserving we
must have that $p_x = -q_y$. Hence $P_2$ must be of the form

\begin{eqnarray} \label{one}
P_2(x, y) = (x^2 + 2b xy + c y^2, -2xy - b y^2).
\end{eqnarray}

From here on we will assume that $b \neq 0$, the case $b=0$ is
almost identical. We can now further simplify Equation
\eqref{one} by conjugating with $(x, y) \rightarrow (x, b^{-1} y)$
to get

\begin{eqnarray}\label{two}
P_2(x,y) = (x^2 + 2xy + c y^2, -2xy -y^2).
\end{eqnarray}

As in \cite{hak1} we blow-up $y = ux$ and write $(x_n, u_n) =
F^n(x,u)$ to get

\begin{align} \label{three}
x_1 = x + (1 + 2u + cu^2) x^2 + O(|x|^3) \\
u_1 = u - (3u + 3u^2 + cu^3) x + O(|x|^2) \label{four}
\end{align}

Recall that in the hyperbolic case the stable manifold is locally
a graph over the attracting direction, and that the mapping is
expanding in the repelling direction. The expansion guarantees
that the local stable manifold is stable under small
perturbations. We will see that the situation is analogous for
volume preserving automorphisms tangent to the identity: for some
$\e>0$ small enough the attracting set is a graph over the set
$\{x \in \CC \mid \max(|x|, |\arg(x) - \pi|) < \e \}$.

For $\e > 0$ define
\begin{align}
W_\e = \{ (x, u) \in \CC^2 \mid \max(|x|, |\arg(x) - \pi|) < \e,
2|u| < |x|\}
\end{align}

The following Lemma follows from Equations \eqref{three} and
\eqref{four}.

\begin{lemma}\label{expansion}
Let $\e>0$ small enough and let $(x,u) , (\tilde{x}, \tilde{u})
\in W_\e$. Suppose that $2|x - \tilde{x}| < |u - \tilde{u}|$. Then
$|u_1 - \tilde{u}_1| > \max(|u - \tilde{u}|, 2|x_1
-\tilde{x}_1|)$.
\end{lemma}

We also have

\begin{lemma}\label{convergence}
Let $\e >0$ be small enough and $(x,u) \in W_\e$. If $(x_n, u_n)
\in W_\e$ for every $n$ then $(x_n, u_n) \rightarrow (0,0)$.
\end{lemma}
\begin{proof} It follows from \eqref{three} that $\{x_n\}$ must
converge to $0$. But then $u_n \rightarrow 0$ by definition of
$W_\e$.
\end{proof}

\begin{proposition}\label{graph}
Let $\e >0$ be small enough and $x \in \CC$ satisfy $|x| < \e$ and
$|\arg(x) - \pi| < \e$. Then there is exactly one $u \in \CC$ such
that $(x_n, u_n) \in W_\e$.
\end{proposition}
\begin{proof}
To show existence, let $U_n = \{ u \in \CC \mid (x_j, u_j) \in
W_\e, j = 1, \cdots , n\}$. It follows from Lemma \ref{expansion}
that $\{U_n\}$ is a nested sequence of non-empty relatively
compact sets, hence the intersection is not empty.

To prove uniqueness, suppose for the purpose of contradiction that
there exist two such points, $u$ and $v$. Then Lemma
\ref{expansion} shows inductively that $|u_n - v_n| > |u-v|$ for
every $n$, but Lemma \ref{convergence} shows that both $u_n$ and
$v_n$ must converge to the origin, so we have a contradiction.
\end{proof}

So we indeed have that the attracting set is locally a graph over
the attracting direction, and we have expansion is the other
direction. Stability follows:

\begin{theorem}\label{tangentstability}
Let $f \in \autidentity$ have a non-degenerate attracting
direction, and let $\Omega$ be the corresponding attracting set.
If $\{f_j\} \subset \autidentity$ converges to $f$, then the
corresponding attracting sets $\Omega_j$ satisfy (locally)
\begin{align}
\mathrm{d}(\Omega_j, \Omega) \rightarrow 0.
\end{align}
\end{theorem}

To prove Theorem \ref{tangentstability}, first assume that for $j$
large enough the maps $f_j$ all have a fixed point at the origin,
and the same non-degenerate direction $(1,0)$. We can do this
because for $j$ large enough this can be assured by conjugating
with an affine map arbitrarily close to the identity.

It follows from Lemmas  \ref{expansion}, \ref{convergence} and
Proposition \ref{graph} that the proof of Lemma \ref{graphs} works
here as well, with expansion now in the $u$-coordinate instead of
the $y$-coordinate.

We obtain the existence of dense attracting sets for volume
preserving automorphisms tangent to the identity.

\begin{theorem}\label{densetangent}
There is a dense $G_\d$-set $\mathcal{S}\subset \autidentity$ such
that each $f\in\mathcal{S}$  has a fixed point tangent to the
identity whose attracting set is dense in $\CC^2$.
\end{theorem}

As the proof is very similar to the proof of Theorem \ref{dense}
we will only outline the differences. The unboundedness of the
attracting sets follows from Theorem 1.10 of \cite{hak2}, which
implies that the attracting set is conformally equivalent to
$\CC$.

To complete Theorem \ref{densetangent} we need to confirm that a
version of Theorem 3.1 from \cite{fo2} holds for automorphisms
tangent to the identity.  Let $\mathcal{V}$  denote the set of
volume preserving automorphisms tangent to the identity at the
origin that have a non-degenerate characteristic direction.  For
each $f\in\mathcal{V}$  we define as in \cite{fo2}:
$$
K_f=\{z\in\CC^2 \mid \{f^n(z)\} \ \mathrm{is \ bounded}\}.
$$
The claim in \cite{fo2} becomes: \emph{There exists a dense $G_\d$
set $\mathcal{V}_1$  in $\mathcal{V}$  such that for every
$f\in\mathcal{V}_1$, the set $K_f$  is an $F_\s$ set with empty
interior}.  We outline the additions that have to be made to prove
the claim.

For $f_0 \in \mathcal{V}$, we define $U_C$ as the interior of the
set of points whose forward orbits are contained in the ball
$B_C$. Note that the set $U_C$  is forward invariant under $f_0$.
Since $f_0$  has a non-degenerate characteristic direction at the
origin, we have a repelling piece of curve (namely, the attracting
set to the origin for $f^{-1}$).  This curve is unbounded (since
it is biholomorphic to $\CC$), so the origin cannot lie not in
$U_C$.

When the existence of the polynomially convex set $X_q\subset U_C$
(as in \cite{fo2}) is established, we may assume that
$X_q\cap\{0\}=\emptyset$.  This follows because $X_q$  is the
orbit of $q$ under the action of a commutative compact Lie group
acting as automorphisms on $U_C$. Hence we may extend the vector
field $\xi$ to be zero on a neighborhood of the origin, and one
can still approximate $\xi$  by divergence free polynomial vector
fields.

When approximating the flow of $\xi$ by the 1-parameter family
$\psi_t$  of volume preserving automorphisms, we may then assume
that $\psi_t$  is tangent to the identity for each $t$, by
composing with the inverses of the derivatives at the origin. The
claim follows just as in \cite{fo2}.

\section{The non-autonomous case}

For many results in complex dynamical systems it makes sense to
ask whether the result also holds in the non-autonomous setting.
Instead of studying the iterations $\{f^n\}$ of a fixed mapping,
one studies the compositions $\{f_n \circ \cdots \circ f_1\}$ for
a sequence of mappings $f_1, f_2, \ldots$. As this gives much more
freedom, it is generally easier to construct (counter-) examples,
but harder to prove that general results still hold in the
autonomous setting.

For a sequence $f_1, f_2, \ldots \in \mathrm{Aut}(\CC^k)$ that all
have a fixed point at the origin one can define the attracting set
as
$$
\O^{\{f_j\}} = \{z \in \CC^k \mid f_n \circ \cdots \circ f_1(z)
\rightarrow 0 \}
$$

The complex structure of such basins has been studied in
\cite{fost}, \cite{pw}, \cite{pe}. The following construction
shows that in the non-autonomous setting one can have an
attracting set contains a neighborhood of the origin for a
sequence of automorphisms that are all tangent to the identity.

\begin{theorem}
There exists a sequence $\{g_j\}_{j=1}^\infty$  of automorphisms
of $\CC^2$ with $g_j(0)=0,\mathrm{d}g_j(0)=\mathrm{Id}$  for all
$j\in\NN$, and such that $\O^{\{g_j\}}=\CC^2$.
\end{theorem}

Instead of $\mathrm{d}g_j(0)=\mathrm{Id}$ we may in fact freely
prescribe the $d$-jets for each of the mappings $\{g_j\}$ and for
any fixed $d$. For any sequence of $d$-jets there exists a
sequence of automorphisms having these $d$-jets and that satisfies
the requirements needed in the proof \cite{fc} \cite{we}.

The same construction works in dimensions higher than $2$.

\begin{proof}

The idea of the proof is simple. We construct an increasing
sequence of subsets $U_1 \subset U_2 \subset \ldots$ of $\CC^2$
whose union is $\CC^2$, and we construct automorphisms $f_1, f_2,
\ldots$ that are all tangent to the identity such that $f_j$ maps
$U_j$ into a (smaller and smaller) neighborhood of the origin.
Then we define $g_1 = f_1$ and $g_j = f_j \circ f_{j-1}^{-1}$ for
$j \ge 2$, and we are done.

For $R >0$ and some small $\e>0$ we define the following subsets
of the complex plane:

\

$\bullet$ \ $\overline\triangle_R = \{ |z| \le R\}$,\

$\bullet$ \ $\Theta_R^\e:=\{re^{i\theta}\in\CC \mid-\e<\theta<\e,
\; 0 \le r\leq R\}$,\

$\bullet$ \ $K_R^\e:=\overline\triangle_R\setminus\Theta_R^\e-\e$
\ (here $-\e$ is a translation to the left),\

$\bullet$ \ $L_R^\e:=\{x+iy\in\CC \mid y=0, \; \e\leq x\leq R\}$.

\

For $\d>0$ small enough we have that the following set is a union
of five disjoint compact sets in $\CC^2$:
$$
N_{R,\d}^{\e}:=\overline
B_\d(0,0)\cup(K_R^\e\times\overline\triangle_R)\cup(L_R^\e\times\{0\})\cup(\{0\}\times
K_R^\e)\cup(\{0\}\times L_R^\e).
$$

We claim that $N_{R,\d}^\e$  is polynomially convex.  To see this,
note first that the set $K_R^\e\times\overline\triangle_R
\cup(\overline\triangle_\d\times\overline\triangle_R)$ is
polynomially convex.  Since $\overline
B_\d(0,0)\subset\overline\triangle_\d\times\overline\triangle_R$
it follows by an application of the Oka-Weil theorem that
$K_R^\e\times\overline\triangle_R\cup\overline B_\d(0,0)$  is
polynomially convex.  In adding the rest of the components, we add
compact sets contained in closed submanifolds of $\CC^2$ (in fact
complex lines), so it suffices to check that the intersection of
$N_{R,\d}^\e$ with these submanifolds is polynomially convex. (to
see that this suffices, first use the defining function, and then
apply the local maximum modulus principle). That the intersection
with the complex lines is polynomially convex follows from Runge's
theorem.

Pick $\r > 0$ small enough such that the balls $\overline
B_\d(0,0)$ and $\overline B_\r(\cdot,\cdot)$  are pairwise
disjoint, and their union is polynomially convex. Then define
automorphisms $\phi_1,...,\phi_5$  of $\CC^2$  such that the
following holds:

\

(i) $\phi_1|_{\overline B_\d(0,0)}=\mathrm{Id}$,\

(ii) $\phi_2(K_R^\e\times\overline\triangle_R)\subset
B_\r(-\e,0)$,\

(iii) $\phi_3(L_R^\e\times\{0\})\subset B_\r(\e,0)$,\

(iv) $\phi_4(\{0\}\times K_R^\e)\subset B_\r(0,-\e)$, \

(v) $\phi_5(\{0\}\times L_R^\e)\subset B_\r(0,\e)$.

\

Choose a small neighborhood $U$  of $N_{R,\d}^\e$  and define a
map $\psi:U\rightarrow\CC^2$  such that $\psi|_{\overline
B_\d(0,0)}=\phi_1,
\psi|_{K_R^\e\times\overline\triangle_R}=\phi_2$, and so forth. By
\cite{fr} the map $\psi$  may be approximated arbitrarily well on
$N_{R,\d}^{\e}$  by an automorphism $f=f_{R,\d}^{\e}$ of $\CC^2$.
By composing with a linear map close to the identity if necessary,
we may assume that $f(0)=0,\mathrm{d}f(0)=\mathrm{Id}$. Note that
if the approximation of $\psi$  is good enough then the image of
$N_{R,\d}^{\e}$  is contained in the ball $B_{2\e}(0,0)$. \

To finish the proof we now choose sequences $R_j,\e_j,\d_j$ for
$j=1,2,...$, such that $\e_j,\d_j,\r_j\searrow 0$  and
$R_j\nearrow\infty$, and such that we can carry out the above
construction for the sets $N_{R_j,\d_j}^{\e_j}$  to get a sequence
of automorphisms $f_j=f_{R_j,\d_j}^{\e_j}$  as above.  Now define
the sequence $g_j$  inductively by $g_1=f_1$ and $g_j=f_j\circ
f_{j-1}^{-1}$  for $j=2,3,...$.

We have that the union of the increasing sequence
$N_{R_j,\d_j}^{\e_j}$ is all of $\CC^2$ and we have
$$
g_j g_{j-1} \cdots
g_1(N_{R_j,\d_j}^{\e_j})=f_j(N_{R_j,\d_j}^{\e_j})\subset
B_{2\e_j}(0,0).
$$

This completes the proof.
\end{proof}
\begin{remark}
Note that the convergence in the above theorem is \emph{pointwise}
- not uniform on compacts.  Specifically we have that points
arbitrarily close to the origin go arbitrarily far out towards
infinity. In light of this example one might ask whether there
exist an attracting set of a non-periodic bounded orbit $p_0, p_1,
\ldots$, such that the attracting set contains a neighborhood of
$p_0$. The arguments in the proof of Theorem \ref{uniform} adapts
to this case however, showing that this is impossible.
\end{remark}

\bibliography{biblio}

\begin{thebibliography}{10}

\bibitem{bf}
G.~T. Buzzard and F.~Forstneri\v{c}.
\newblock An interpolation theorem for holomorphic automorphisms of {${\bf
  C}\sp n$}.
\newblock {\em J. Geom. Anal.}, 10(1):101--108, 2000.

\bibitem{fost}
J.~E. Forn{\ae}ss and B.~Stens{\o}nes.
\newblock Stable manifolds of holomorphic hyperbolic maps.
\newblock {\em Internat. J. Math.}, 15(8):749--758, 2004.

\bibitem{fo2}
J.E. Forn{\ae}ss and N.~Sibony.
\newblock The closing lemma for holomorphic maps.
\newblock {\em Ergod.Th.and Dynam.Sys.}, 17, 1997.

\bibitem{fc}
F.~Forstneri\v{c}.
\newblock Interpolation by holomorphic automorphisms and embeddings in
  $\mathbb{C}^n$.
\newblock {\em J. of Geom. Anal.}, 9, 1999.

\bibitem{fr}
F.~Forstneri\v{c} and J.-P. Rosay.
\newblock Approximation of biholomorphic mappings by automorphisms of
  $\mathbb{C}^n$.
\newblock {\em Invent. Math.}, 112:323--349, 1993.

\bibitem{fr2}
F.~Forstneri\v{c} and J.-P. Rosay.
\newblock Erratum, approximation of biholomorphic mappings by automorphisms of
  $\mathbb{C}^n$.
\newblock {\em Invent. Math.}, 118:573--574, 1994.

\bibitem{hak2}
M.~Hakim.
\newblock Transformations tangent to the identity.
\newblock 1997.
\newblock Stable pieces of manifolds.

\bibitem{hak1}
M.~Hakim.
\newblock Analytic transformations of {$(\bold C\sp p,0)$} tangent to the
  identity.
\newblock {\em Duke Math. J.}, 92(2):403--428, 1998.

\bibitem{kh}
A.~Katok and B.~Hasselblatt.
\newblock {\em Introduction to the modern theory of dynamical systems},
  volume~54 of {\em Encyclopedia of Mathematics and its Applications}.
\newblock Cambridge University Press, Cambridge, 1995.

\bibitem{mo}
S.~Morosawa, Y.~Nishimura, M.~Taniguchi, and T.~Ueda.
\newblock {\em Holomorphic dynamics}, volume~66 of {\em Cambridge Studies in
  Advanced Mathematics}.
\newblock Cambridge University Press, 2000.

\bibitem{pe}
H.~Peters.
\newblock Perturbed basins of attraction.
\newblock {\em Math. Ann.}, 2006.

\bibitem{pw}
H.~Peters and E.~F. Wold.
\newblock Non-autonomous basins of attraction and their boundaries.
\newblock {\em J. Geom. Anal.}, 15(1):123--136, 2005.

\bibitem{rr}
J.P. Rosay and W.~Rudin.
\newblock Holomorphic maps from $\mathbb{C}^n$ to $\mathbb{C}^n$.
\newblock {\em Trans. Amer. Mat. Soc.}, 310:47--86, 1998.

\bibitem{ue}
T.~Ueda.
\newblock Local structure of analytic transformations of two complex variables.
  {I}.
\newblock {\em J. Math. Kyoto Univ.}, 26(2):233--261, 1986.

\bibitem{we}
B.~Weickert.
\newblock {\em Automorphisms of $\mathbb{C}^n$.}
\newblock PhD thesis, University of Michigan, 1997.

\end{thebibliography}
\end{document}